\documentclass[MSNbibl,number,citesort,seceqn,dvips]{arxbj}

\aid{0}
\volume{20}
\issue{1}
\pubyear{2014}
\firstpage{207}
\lastpage{230}
\doi{10.3150/12-BEJ482} 

\makeatletter

\newcommand{\RMO}{\mathrm{O}}

\newcommand{\RMe}{\mathrm{e}}

\newcommand{\mrmd}{\mathrm{d}}

\newtheorem{theorem}{Theorem}[section]
\newtheorem{cor}{Corollary}[section]
\newtheorem{lem}{Lemma}[section]

\newremark{rem}{Remark}[section]

\newproclaim{assumption}{Assumption}[section]

\newremark{example}{Example}

\makeatother

\begin{document}
\begin{frontmatter}

\title{Uniform convergence rates for a class of martingales with
application in non-linear cointegrating regression}
\runtitle{Uniform convergence rates for a class of martingales}

\begin{aug}
\author{\fnms{Qiying} \snm{Wang}\corref{}\thanksref{e1}\ead[label=e1,mark]{qiying.wang@sydney.edu.au}} \and
\author{\fnms{Nigel} \snm{Chan}\thanksref{e2}\ead[label=e2,mark]{chanhiungai@gmail.com}}
\runauthor{Q. Wang and N. Chan} 
\address{School of Mathematics and Statistics, The University of
Sydney, NSW 2006,
Australia.\\ \printead{e1,e2}}
\end{aug}

\received{\smonth{10} \syear{2011}}
\revised{\smonth{9} \syear{2012}}

%
\begin{abstract}
For a class of martingales, this paper provides a framework on the
uniform consistency with broad applicability. The main condition
imposed is only related to the conditional variance of the martingale,
which holds true for stationary mixing time series, stationary iterated
random function, Harris recurrent Markov chains and $I(1)$ processes
with innovations being a linear process. Using the established results,
this paper investigates the uniform convergence of the Nadaraya--Watson
estimator in a non-linear cointegrating regression model. Our results
not only provide sharp convergence rate, but also the optimal range
for the uniform convergence to be held. This paper also considers the
uniform upper and lower bound estimates for a functional of Harris
recurrent Markov chain, which are of independent interests.
\end{abstract}

%

\end{frontmatter}

\section{Introduction}
Let $(u_k, x_k)$ with $x_k=(x_{k1},\ldots, x_{kd}), d\ge1$, be a sequence
of random vectors. A common functional of interests $S_n(x)$ of $(u_k,
x_k)$ is defined by
%
%
\begin{equation}\label{in1}
S_n(x) =\sum_{k=1}^n
u_k f\bigl[(x_k+x)/h\bigr],\qquad x\in R^d,
\end{equation}
where $h=h_n\to0$ is a certain sequence of positive constants and
$f(x)$ is a real function on~$R^d$. Such functionals arise in
non-parametric estimation problems, where
$f$ may be a kernel function $K$ or a squared kernel function $K^{2}$
and the
sequence $h$ is the bandwidth used in the non-parametric regression.

The uniform convergence of $S_n(x)$ in the situation that the $(u_k,
x_k)$ satisfy certain stationary conditions was studied in many
articles. Liero \cite{liero}, Peligrad \cite{peligrad} and Nze and
Doukhan \cite{angonze} considered the uniform convergence over a fixed
compact set, while Masry \cite{masry}, Bosq \cite{bosq} and Fan and
Yao \cite{fan} gave uniform results over an unbounded set. These work
mainly focus on random sequence $x_t$ which satisfies different types
of mixing conditions. Investigating a more general framework, Andrews
\cite{andrew} gave result on kernel estimate when the data sequence is
near-epoch dependent on another underlying mixing sequence. More
recently, Hansen \cite{hansen} provided a set of general uniform
consistency results, allowing for stationary strong mixing multivariate
data with infinite support, kernels with unbounded support and general
bandwidth sequences. Kristensen \cite{kristensen} further extended
Hansen's results to the heterogenous dependent case under $\alpha
$-mixing condition. Also see Wu, Huang and Huang \cite{wu1} for kernel
estimation in general time series settings.

In comparison to the extensive results where the $x_k$ comes from a
stationary time series data, there is little investigation on the the
uniform convergence of $S_n(x)$ for the $x_k$ being a non-stationary
time series. In this regard, Gao, Li and Tj{\o}stheim \cite{gao3}
derived strong and weak consistency results for the case where the
$x_k$ is a null-recurrent Markov chain. Wang and Wang \cite{wangwang}
worked with partial sum processes of the\vspace*{1pt} type $x_k=\sum_{j=1}^k\xi_j$
where $\xi_j$ is a general linear process. While the rate of
convergence in Gao, Li and Tj{\o}stheim \cite{gao3} is sharp, they
impose the independence between $u_k$ and $x_k$. Using a quite
different method, Wang and Wang \cite{wangwang} allowed for the
endogeneity between $u_k$ and $x_k$, but their results hold only for
the $x$ being in a fixed compact set.

The aim of this paper is to present a general uniform consistency
result for $S_n(x)$ with broad applicability. As a framework, our
assumption on the $x_t$ is only related to the conditional variance of
the martingale, that is, $ \sum_{t=1}^n f^2[(x_t+x)/h]$. See
Assumption \ref{assumption2.3} in Section \ref{sec2}. This of course is
a ``high
level'' condition, but it in fact is quite natural which holds true for
many interesting and important examples, including stationary mixing
time series, stationary iterated random function and Harris recurrent
Markov chain. See Sections \ref{sec22} and \ref{sec23} for the
identification of Assumption~\ref{assumption2.3}. This condition also
holds true for
$I(1)$ processes with innovations being a linear process, but the
identification is complicated and requires quite different techniques.
We will report related work in a separate paper. By using the
established result, we investigate the uniform convergence of the
Nadaraya--Watson estimator in a non-linear cointegrating regression
model. It confirms that the uniform asymptotics in Wang and Wang \cite
{wangwang} can be extended to a unbounded set and the independence
between the $u_t$ and $x_t$ in Gao, Li and Tj{\o}stheim \cite{gao3}
can be removed. More importantly, our result not only provides sharp
convergence rate, but also the optimal range for the uniform
convergence to be held. It should be mentioned that our work on the
uniform upper and lower bound estimation for a functional of Harris
recurrent Markov chain is of independent interests.

This paper is organized as follows. Our main results are presented in
next section, which includes the establishment of a framework on the
uniform convergence for a class of martingale and uniform upper and
lower bound estimation for a functional of Harris recurrent Markov
chain. An application of the main results in non-linear cointegrating
regression is given in Section \ref{sec3}. All proofs are postponed to
Section \ref{sec4}. Throughout the paper, we denote constants by $C,
C_1, C_2,\ldots$ which may be different at each appearance. We also use
the notation $\|x\|=\max_{1\le i\le d}|x_i|$.

\section{Main results} \label{sec2}

\subsection{Uniform convergence for a class of martingales} \label{sec21}
We make use of the following assumptions in the development of uniform
convergence for the $S_n(x)$ defined by (\ref{in1}).
Recall $x_k=(x_{k1},\ldots, x_{kd})$ where $ d\ge1$ is an integer.
%
%
\begin{assumption}\label{assumption2.1}
$\{u_t, {\mathcal F}_t\}_{t\ge1}$ is a
martingale difference, where ${\mathcal F}_t=\sigma(x_1,\ldots,
x_{t+1},\break
u_1,\ldots,u_t)$, satisfying $ \sup_{t\ge1}E(|u_t|^{2p}\mid{\mathcal
F}_{t-1})<\infty$, a.s., for some $p\ge1$ specified in Assumption
\ref{assumption2.4} below.
\end{assumption}
%
%
\begin{assumption}\label{assumption2.2}
$f(x)$ is a real function on $R^d$ satisfying
$\sup_{x\in R^d} |f(x)|<\infty$ and
$|f(x)-f(y)| \le C \|x-y\| $ for all $x, y\in R^d$ and some constant $C>0$.
\end{assumption}
%
%
\begin{assumption}\label{assumption2.3}
There exist positive constant sequences
$c_n\uparrow\infty$ and $b_n$ with $b_n=\RMO(n^k)$ for some $k>0$ such that
%
%
\begin{equation}\label{m1}
\sup_{\|x\|\le b_n} \sum_{t=1}^n
f^2\bigl[(x_t+x)/h\bigr] =\RMO_P(c_n).
\end{equation}
\end{assumption}
%
%
\begin{assumption}\label{assumption2.4}
$h\to0, nh\to\infty$ and $n c_n^{-p} \log^{p-1}n=\RMO(1)$, where $c_n$ is
defined as in Assumption \ref{assumption2.3} and $p$ is defined as in
Assumption \ref{assumption2.1}.
\end{assumption}

We remark that Assumption \ref{assumption2.1} ensures that $\{S_n(x),
{\mathcal F}_n\}_{n\ge1}$ is a martingale
for each fixed $x$ and is quite weak. Clearly, Assumption \ref
{assumption2.1} is
satisfied if $u_t$ is a sequence of i.i.d. random variables, which is
independent of $x_1,\ldots, x_t$, with $Eu_1=0$ and $E|u_1|^{2p}<\infty$.
The Lipschitz condition used in Assumption \ref{assumption2.2} is
standard in the
investigation of uniform consistency, where we do not require the
$f(x)$ to have finite compact support.
Assumption \ref{assumption2.3} is a ``high level'' condition for the
$x_k$. We use it
here to provide a framework. In Sections \ref{sec22} and \ref
{sec23}, we will show that this condition is in fact quite natural
which holds true by many interesting and important examples. Assumption
\ref{assumption2.4} provides the connections among the moment condition
required in
Assumption~\ref{assumption2.1}, the condition (\ref{m1}) and the
bandwidth $h$. In
many applications, we have
$c_n=n^{\alpha} h^{d} l(n)$, where $ 0<\alpha\le1$ and
$l(n)$ is a slowly varying function at infinite. See Section \ref
{sec23} and Examples \ref{example1}--\ref{example3} in Section \ref
{sec22}. In the typical
situation that $c_n=n^{\alpha} h^{d} l(n)$, if there exists a
\mbox{$0<\varepsilon_0<\alpha$} such that $n^{\alpha-\varepsilon_0}h^{d}\to
\infty$, the $p$
required in Assumption \ref{assumption2.1} can be specified to
$p=[1/\varepsilon_0]+1$.

We have the following main result.
%
%
\begin{theorem} \label{th4} Under Assumptions \ref{assumption2.1}--\ref
{assumption2.4}, we have
%
%
\begin{equation}\label{m10}
\sup_{\|x\|\le b_n} \Biggl| \sum_{t=1}^{n}
u_t f\bigl[(x_t+x)/h\bigr] \Biggr| =\RMO_P
\bigl[(c_n \log n)^{1/2} \bigr].
\end{equation}
If (\ref{m1}) is replaced by
%
%
\begin{equation}\label{m1a9}
\sup_{\|x\|\le b_n} \sum_{t=1}^n
f^2\bigl[(x_t+x)/h\bigr] = \RMO(c_n),\qquad\mbox{a.s.},
\end{equation}
the result (\ref{m10}) can be strengthened to
%
%
\begin{equation}\label{m11}
\sup_{\|x\|\le b_n} \Biggl| \sum_{t=1}^{n}
u_t f\bigl[(x_t+x)/h\bigr] \Biggr| =\RMO \bigl[(c_n
\log n)^{1/2} \bigr],\qquad\mbox{a.s.}
\end{equation}
\end{theorem}

Theorem \ref{th4} can be extended to uniform convergence for the
$S_n(x)=\sum_{t=1}^{n} u_t f[(x_t+x)/h]$ over unrestricted space
$R^d$. This requires additional condition on the $x_k$ and the tail
decay for the function $f(x)$.
%
%
\begin{theorem} \label{th4a} In addition to Assumptions \ref
{assumption2.1}--\ref{assumption2.4}, $n
\sup_{\|x\|> b_n/2} |f(x/h)| =\break \RMO [(c_n \log n)^{1/2} ]$ and there
exists a $k_0>0$ such that
%
%
\begin{equation}\label{p09}
b_n^{-k_0} \sum_{t=1}^n
E\|x_t\|^{k_0}= \RMO \bigl[(c_n \log
n)^{1/2} \bigr].
\end{equation}
Then,
%
%
\begin{equation}\label{m10a}
\sup_{x\in R^d} \Biggl| \sum_{t=1}^{n}
u_t f\bigl[(x_t+x)/h\bigr] \Biggr| =\RMO_P
\bigl[(c_n \log n)^{1/2} \bigr].
\end{equation}
Similarly, if
(\ref{m1}) is replaced by (\ref{m1a9}) and (\ref{p09}) is replaced by
%
%
\begin{equation}\label{p0}
b_n^{-k_0} \sum_{t=1}^n
\|x_t\|^{k_0}= \RMO \bigl[(c_n \log
n)^{1/2} \bigr],\qquad\mbox{a.s.},
\end{equation}
then
%
%
\begin{equation}\label{m11a}
\sup_{x\in R^d} \Biggl| \sum_{t=1}^{n}
u_t f\bigl[(x_t+x)/h\bigr] \Biggr| =\RMO \bigl[(c_n
\log n)^{1/2} \bigr],\qquad\mbox{a.s.}
\end{equation}
\end{theorem}
%
%
\begin{rem} \label{rem21} Theorems \ref{th4}--\ref{th4a} allow for
the $x_t$ to be a stationary or non-stationary time series. See
Examples \ref{example1}--\ref{example3} and Section \ref{sec23} below.
More examples on
non-stationary time series will be reported in a separate paper. The
rates of convergence in both theorems are sharp. For instance, in the
well-known stationary situation such as those appeared in Examples
\ref{example1}--\ref{example3}, the $c_n$ can be chosen as $c_n=nh$.
Hence, when there are enough
moment conditions on the $u_t$ (i.e., $p$ is large enough), we obtain
the optimal rate $n^{2/5}\log^{3/5} n$, by taking $h\sim(\log
n/n)^{1/5}$. In non-stationary situation, the rate of convergence is
different. In particular we have $c_n=\sqrt nh$ for the $x_t$ to be a
random walk given in Corollary \ref{th3}. The reason behind this fact
is that the amount of time spent by the random walk around any
particular point is of order $\sqrt n$ rather than $n$ for a stationary
time series. For more explanation in this regard, we refer to
Wang and Phillips \cite{wangphillips1,wangphillips2}.
\end{rem}

\subsection{\texorpdfstring{Identifications of Assumption \protect\ref{assumption2.3}}
{Identifications of Assumption 2.3}}
\label{sec22} This
section provides several stationary time series examples which satisfy
Assumption \ref{assumption2.3}. Examples \ref{example1} and
\ref{example2} come from Wu, Huang and Huang \cite{wu1}, where more
general settings on the $x_t$ are established. Example \ref{example3}
discusses a strongly mixing time series. This example comes from Hansen
\cite{hansen}. By making use of other related works such as Peligrad
\cite{peligrad}, Nze and Doukhan~\cite{angonze}, Masry \cite{masry},
Bosq \cite{bosq} and Andrews \cite{andrew}, similar results can be
established for other mixing time series like $\rho$-mixing and
near-epoch-dependent time series. In these examples, we only consider
the situation that $d=1$. The extension to $d> 1$ is straightforward
and hence the details are omitted. Throughout Examples
\ref{example1}--\ref{example3}, we use the notation $f^2_h(x)=
h^{-1}f^2(x/h)$.

Example on the Harris recurrent Markov chains, which allows for
stationary (positive recurrent) or non-stationary (null recurrent)
series, is given in Section \ref{sec23}. In the section, we also
consider the uniform lower bound, which is of independent interests.
More examples on $I(1)$ processes with innovations being linear
processes will be reported in a separate paper.
%
%
\begin{example}\label{example1}
Let $\{x_{t}\}_{t\geq0}$ be a linear process defined by
\[
x_{t}=\sum_{k=0}^{\infty}
\phi_{k} \varepsilon_{t-k},
\]
where $\{\varepsilon_{j} \}_{j\in Z}$ is a sequence of i.i.d.
random variables with $E\varepsilon_{0}^{2}<\infty$ and a density
$p_\varepsilon
$ satisfying $\sup_x|p_\varepsilon^{(r)}(x)|<\infty$ and
\[
\int_{ R} \bigl|p_\varepsilon^{(r)}(x)\bigr|^2
\,\mrmd x< \infty,\qquad r=0, 1, 2,
\]
where $p_\varepsilon^{(r)}(x)$ denotes the $r$-order derivative of
$p_\varepsilon
(x)$. Suppose that $\sum_{k=0}^\infty|\phi_k| < \infty$ and $\phi
\equiv\sum_{k=0}^\infty\phi_k \ne0$, and in addition Assumption
\ref{assumption2.2}, $f(x)$ has a compact support. It follows from
Section 4.1 of Wu,
Huang and Huang \cite{wu1} that,
for any $h \rightarrow0$ and $nh \log^{-1}n \rightarrow\infty$,
%
%
\begin{equation}
\sup_{x \in R} \Biggl| \frac{1}{n} \sum_{t =1}^{n}
\bigl[f_h^2(x_t + x) - Ef_h^2(x_t
+ x) \bigr] \Biggr| = \RMO \biggl[\sqrt{\frac{\log n}{nh}} + n^{-1/2} l(n) \biggr]
,\qquad\mbox{a.s.},
\end{equation}
where $l(n)$ is a slowly varying function.
Note that $x_t$ is stationary process with a bounded density $g(x)$
under the given conditions on $\varepsilon_k$. Simple calculations show that
%
%
\begin{equation}
\sup_{x \in R} \sum_{t =1}^{n}
f^2\bigl[(x_t + x)/h\bigr] = \RMO_P(nh),
\end{equation}
that is, $x_t$ satisfies Assumption \ref{assumption2.3}.\vadjust{\goodbreak}
\end{example}
%
%
\begin{example}\label{example2}
Consider the non-linear time series of the following form
\begingroup
\abovedisplayskip=7pt
\belowdisplayskip=7pt
\[
x_k = R(x_{k-1}, \varepsilon_k),
\]
where $R$ is a bivariate measurable function and $\varepsilon_k$ are i.i.d.
innovations. This is the iterated random function framework that
encompasses a lot of popular non-linear time series models. For example,
if $R(x, \varepsilon) = a_1 x I(x < \tau) + a_2 x I(x \ge\tau) +
\varepsilon$, it is the threshold autoregressive (TAR) model (see Tong
\cite
{tong}). If $R(x, \varepsilon) = \varepsilon\sqrt{a_1^2 + a_2^2 x}$,
then it is
autoregressive model with conditional heteroscedasticity (ARCH) model.
Other non-linear time series models, including random coefficient
model, bilinear autoregressive model and exponential autoregressive
model can be fitted in this framework similarly. See Wu and Shao \cite
{wu2} for details.

In order to identify Assumption \ref{assumption2.3}, we need some
regularity conditions
on the initial distribution of $x_0$ and the function $R(x, \varepsilon
)$. Define
%
%
\begin{equation}
L_{\varepsilon} = \sup_{x \ne x'} \frac{ | R(x, \varepsilon) - R(x',
\varepsilon) | } { |x
- x'| }.
\end{equation}
Denote by $g(x\mid x_0)$ the conditional density of $x_1$ at $x$ given
$x_0$. Further let $g'(y\mid x) = \partial g(y\mid x) / \partial y$ and
%
%
\begin{equation}
I(x) = \biggl[ \int_{{R}} \biggl| \frac{\partial}{\partial x} g(y\mid x)
\biggr|^2 \,\mrmd y \biggr]^{1/2} \quad\mbox{and}\quad J(x) = \biggl[ \int
_{{R}} \biggl| \frac{\partial}{\partial x} g'(y\mid x)
\biggr|^2 \,\mrmd y \biggr]^{1/2},
\end{equation}
$I(x)$ and $J(x)$ can be interpreted as a prediction sensitivity
measure. These quantities measure the change in 1-step predictive
distribution of $x_1$ with respect to change in initial value $x_0$.
Suppose that:
\begin{longlist}[(iii)]
\item[(i)] there exist $\alpha$ and $z_0$ such that
\[
E \bigl( |L_{\varepsilon_0}|^{\alpha} + \bigl| R(z_0,
\varepsilon_0)\bigr|^\alpha\bigr)<\infty,\qquad E \bigl[
\log(L_{\varepsilon_0}) \bigr] < 0 \quad\mbox{and}\quad E L_{\varepsilon
_0}^2 < 1;
\]
\item[(ii)]
$\sup_x [I(x) + J(x)] < \infty$;
\item[(iii)] in addition to Assumption \ref{assumption2.2}, $f(x)$ has
a compact support.
\end{longlist}
It follows from Section 4.2 of Wu, Huang and Huang \cite{wu1} that,
for any $h \rightarrow0$ and $nh \log^{-1}n \rightarrow\infty$
%
%
\begin{equation}
\sup_{x \in R} \Biggl| \frac{1}{n} \sum_{t =1}^{n}
\bigl[f_h^2(x_t + x) - Ef_h^2(x_t
+ x) \bigr] \Biggr| = \RMO \biggl[\sqrt{\frac{\log n}{nh}} + n^{-1/2} l(n) \biggr]
,\qquad\mbox{a.s.},
\end{equation}
where $l(n)$ is a slowly varying function.
Note that $x_t$ has a unique and stationary distribution under the
given condition (i) and (ii). See Diaconis and Freedman \cite
{diaconis}, for instance. Simple calculations show that
%
%
\begin{equation}
\sup_{x \in R} \sum_{t =1}^{n}
f^2\bigl[(x_t + x)/h\bigr] = \RMO_P(nh),
\end{equation}
that is, $x_t$ satisfies Assumption \ref{assumption2.3}.\vadjust{\goodbreak}
\end{example}\endgroup
%
%
\begin{example}\label{example3}
Let $\{x_k\}_{k\ge0}$ be a strictly stationary time
series with density $g(x)$. Suppose that:
\begin{longlist}[(iii)]
\item[(i)] $x_t$ is strongly mixing with mixing coefficients $\alpha(m)$
that satisfy $\alpha(m) \le A m^{-\beta}$ where $\beta> 2 $ and $A <
\infty$;
\item[(ii)] $\sup_x|x|^q g(x) < \infty$ for some $q\ge1$ satisfying
$\beta>2+1/q$ and there is some $j^* < \infty$ such that for all $j
\ge j^*$, $\sup_{x, y} g_j(x, y) < \infty$ where $g_j(x, y)$ is the
joint density of $\{x_0, x_j\}$;
\item[(iii)] in addition to Assumption \ref{assumption2.2}, $f(x)$ has
a compact support.
\end{longlist}
It follows from Theorem 4 (with $Y_i=1$) of Hansen \cite{hansen} that,
for any $h \rightarrow0$ and $n^{\theta}h \log^{-1}n \rightarrow
\infty$
with $\theta=\beta-2-1/q$,
%
%
\begin{equation}\label{exa3}
\sup_{x \in R} \Biggl| \frac{1}{n} \sum_{t =1}^{n}
\bigl[f_h^2(x_t + x) - Ef_h^2(x_t
+ x) \bigr] \Biggr| = \RMO_P \biggl[\sqrt{\frac{\log
n}{nh}} \biggr].
\end{equation}
If in addition $E|x_0|^{2q}<\infty$, the result (\ref{exa3}) can be
strengthened to almost surely convergence. Simple calculations show that
%
%
\begin{equation}
\sup_{x \in R} \sum_{t =1}^{n}
f^2\bigl[(x_t + x)/h\bigr] = \RMO_P(nh),
\end{equation}
that is, $x_t$ satisfies Assumption \ref{assumption2.3}.
\end{example}

\subsection{Uniform bounds for functionals of Harris recurrent Markov
chain} \label{sec23}

Let $\{x_k\}_{k\ge0}$ be a Harris recurrent Markov chain with state
space $(E, \mathcal{E})$,
transition probability $P(x, A)$ and invariant measure $\pi$. We
denote $P_\mu$ for the Markovian probability
with the initial distribution $\mu$, $E_\mu$ for correspondent
expectation and $P^k(x, A)$
for the $k$-step transition of $\{x_k\}_{k\ge0}$. A subset $D$ of $E$
with $0<\pi(D)<\infty$ is called $D$-set of $\{x_k\}_{k\ge0}$ if for
any $A\in\mathcal{E}^+$,
\[
\sup_{x\in E} E_x \Biggl(\sum_{k=1}^{\tau_A}I_D(x_k)
\Biggr)<\infty,
\]
where $ \mathcal{E}^+=\{A\in\mathcal{E}\dvtx \pi(A)>0\}$ and $\tau
_A=\inf\{n\ge1\dvtx x_n\in A\}$. As is well-known,
$D$-sets not only exist, but generate the entire sigma $\mathcal{E}$, and
for any $D$-sets $C, D$ and any probability measure $\nu, \mu$ on
$(E, \mathcal{E})$,
%
%
\begin{equation}\label{13}
\lim_{n\to\infty}\sum_{k=1}^n\nu
P^k(C)\Big/\sum_{k=1}^n\mu
P^k(D) =\frac{\pi(C)}{\pi(D)},
\end{equation}
where $\nu P^k(D) =\int_{-\infty}^{\infty} P^k(x, D)\nu(\mrmd x)$. See
Nummelin \cite{nummelin}, for instance.

Let a $D$-set $D$ and a probability measure $\nu$ on $(E, \mathcal
{E})$ be fixed. Define
\[
a(t) = \pi^{-1}(D)\sum_{k=1}^{[t]}
\nu P^k(D),\qquad t\ge0.
\]
By recurrence, $a(t)\to\infty$.
By virtue of (\ref{13}), the asymptotic order of $a(t)$
depends only on $\{x_k\}_{k\ge0}$. As in Chen \cite{chen}, a Harris
recurrent Markov chain $\{x_k\}_{k\ge0}$ is called
$\beta$-regular if
%
%
\begin{equation}\label{d1}
\lim_{\lambda\to\infty} a(\lambda t)/a(\lambda) =t^\beta\qquad \forall t>0,
\end{equation}
where $0< \beta\le1$. It is interesting to notice that, under the
condition (\ref{d1}), the function $a(t)$ is regularly varying at
infinity, that is, there exists a slowly varying function $l(x)$ such that
$a(t)\sim t^{\beta} l(t)$. This implies that the definition of
$\beta$-regular Harris recurrent Markov chain is similar to that of
$\beta$-null recurrent given in Karlsen and Tj{\o}stheim \cite
{karlsen1} and Gao, Li and Tj{\o}stheim \cite{gao3},
but it is more natural and simple.

The following theorem provides uniform upper and lower bounds for a
functional of $x_t$. The upper bound implies that $x_t$ satisfies
Assumption \ref{assumption2.3}, allowing for the $x_t$ being stationary
($\beta=1$,
positive recurrent Markov chain) and non-stationary ($0<\beta<1$, null
recurrent Markov chain). The lower bound plays a key role in the
investigation of the uniform consistency for the kernel estimator in a
non-linear cointegrating regression, and hence is of independent
interests. See Section \ref{sec3} for more details. Both upper and
lower bounds are optimal, which is detailed in Remarks \ref{rem22}
and \ref{rem23}.
%
%
\begin{theorem} \label{th1} Suppose that:
\begin{longlist}[(ii)]
\item[(i)] $\{x_k\}_{k\ge0}$ is a $\beta$-regular Harris recurrent
Markov chain,
where the invariant measure $\pi$ has a bounded density function
$p(s)$ on $R$;
\item[(ii)] in addition to Assumption \ref{assumption2.2}, $\int
_{-\infty}^{\infty
}|f(x)|\,\mrmd x<\infty$.
\end{longlist}
Then, for any $h>0$ satisfying $n^{-\varepsilon_0}a(n) h\to\infty$ for
some $\varepsilon_0>0$, we have
%
%
\begin{equation}\label{12}
\sup_{|x|\le n^m} \sum_{k=1}^n
f^2\bigl[(x_k+x)/h\bigr] =\RMO_P\bigl[a(n) h
\bigr],
\end{equation}
where $m$ can be any finite integer.

For a given sequence of constants $b_n>0$, if there exists a constant
$C_0>0$ such that, uniformly for $n$ large enough,
%
%
\begin{equation}\label{wond}
\inf_{|x|\le b_n+1} \sum_{k=1}^n E
f^2\bigl[(x_k+x)/h\bigr]\ge a(n) h/C_0,
\end{equation}
then, for any $h>0$ satisfying $n^{-\varepsilon_0}a(n) h\to\infty$ for
some $\varepsilon_0>0$, we have
%
%
\begin{equation}\label{12a}
\Biggl\{\inf_{|x|\le b_n} \sum_{k=1}^n
f^2\bigl[(x_k+x)/h\bigr] \Biggr\}^{-1}
=\RMO_P \bigl\{\bigl[a(n) h\bigr]^{-1} \bigr\}.
\end{equation}
\end{theorem}
%
%
\begin{rem} \label{rem22} The result (\ref{12a}) implies that, for
any $0<\eta<1$, there exists a constant $C_\eta>0$ such that
%
%
\begin{equation}\label{15}
P \Biggl(\inf_{|x|\le b_n} \sum_{k=1}^n
f^2\bigl[(x_k+x)/h\bigr] \ge a(n) h/C_\eta
\Biggr) \ge1-\eta.
\end{equation}
This makes both bounds on (\ref{12}) and (\ref{12a}) are optimal.
On the other hand, since the result (\ref{15}) implies that
\[
E \inf_{|x|\le b_n} \sum_{k=1}^n
f^2\bigl[(x_k+x)/h\bigr]\ge a(n) h (1-\eta
)/C_{\eta}
\]
for any $0<\eta<1$, the condition (\ref{wond}) is close to minimal.
\end{rem}

Note that random walk is a $1/2$-regular Harris recurrent Markov chain.
The following corollary on a random walk shows the range $|x|\le b_n$
can be taken to be optimal as well.
%
%
\begin{cor} \label{th3} Let $\{\varepsilon_{j}, 1 \le j\le n \}$ be a
sequence of i.i.d.
random variables with $E\varepsilon_{0}=0$, $E\varepsilon_{0}^{2}=1$
and the
characteristic function $\varphi(t)$ of $\varepsilon_{0}$ satisfying
$\int_{-\infty
}^{\infty}|\varphi(t)|\,\mrmd t<\infty$. Write $x_t=\sum_{j=1}^t\varepsilon_j,
t\ge1$. If in addition to Assumption \ref{assumption2.2}, $\int_{-\infty
}^{\infty
}|f(x)|\,\mrmd x<\infty$, then, for $h> 0$ and $n^{1/2-\varepsilon_0} h\to
\infty
$ where $0<\varepsilon_0<1/2$, we have
%
%
\begin{equation}\label{12a1}
\sup_{|x|\le n^m} \sum_{k=1}^n
f^2\bigl[(x_k+x)/h\bigr] = \RMO_P(\sqrt n h)
\end{equation}
for any integer $m>0$, and
%
%
\begin{equation}\label{m15}
\Biggl\{\inf_{|x|\le\tau_n\sqrt n} \sum_{k=1}^n
f^2\bigl[(x_k+x)/h\bigr] \Biggr\}^{-1} =
\RMO_P \bigl\{(\sqrt n h)^{-1} \bigr\}
\end{equation}
for any $0<\tau_n\to0$.
\end{cor}
%
%
\begin{rem} \label{rem23} For a random walk $x_t$ defined as in
Corollary \ref{th3}, it was shown in Wang and Phillips \cite
{wangphillips1} that
%
%
\begin{equation}\label{ez1}
\frac{1}{\sqrt nh} \sum_{t=1}^{n}
f^2 \bigl[(x_t + y_n )/h\bigr]
\to_D\int f^2(s) \,\mrmd s L_W(1,y),
\end{equation}
where $L_{W}(1,y)$ is a local time of a Brownian motion $W_t$, and
$y=0$ if $y_n/\sqrt n\to0$ and $y=y_0$ if $y_n/\sqrt n\to y_0$. Since
$P(L_W(1,y)=0)>0$ for any $y\not=0$,
it follows from (\ref{ez1}) that the range $\inf_{|x| \le\tau_n
\sqrt{n}}$ in (\ref{m15}) cannot be extended to $\inf_{|x|\le
d\sqrt n}$ for any $d>0$.
\end{rem}
%
%
\begin{rem} \label{rem24} As in Examples \ref{example1}--\ref
{example3}, we may obtain a better
result if $\{x_t\}_{t\ge0}$ is stationary (positive null recurrent)
and satisfies certain other restrictive conditions. Indeed, Kristensen
\cite{kristensen} provided such a result.\vadjust{\goodbreak}

Let $\{x_n\}_{n\ge0}$ be a time-homogeneous, geometrically ergodic
Markov chain. Denote the 1-step transition probability by $p(y\mid x)$,
such that $P(x_{i+1} \in A \mid x_i) = \int_A p(y\mid x) \,\mrmd y$. Also denote
the $i$-step transition probability by $p_i(y\mid x)$, such that
$p_i(y\mid x)
= \int_{ R} p(y\mid z)p_{i-1}(z\mid x) \,\mrmd z$.
Since $x_t$ is geometrically ergodic, it has a density $g(x)$. Further
suppose that:
\begin{longlist}[(iii)]
\item[(i)] (strong Doeblin condition) there exists $s \ge1$ and $\rho
\in(0, 1)$ such that for all $y \in R$,
%
%
\begin{equation}
p_s(y \mid x) \ge\rho g(y);
\end{equation}
\item[(ii)] $\partial^r p(y\mid x) / \partial y^r$ exists and is
uniformly continuous for all $x$, for some $r \ge1$,
\item[(iii)] $\sup_y [ g(y) + |y|^q g(y)]<\infty$ for some $q \ge1$,
\item[(iv)] in addition to Assumption \ref{assumption2.2}, $f(x)$ has a
compact support.
\end{longlist}
It follows from Kristensen \cite{kristensen} that, for any $h\to0$
and $nh\to\infty$,
%
%
\begin{equation}\label{ex3}
\sup_{x \in R} \Biggl| \frac{1}{nh} \sum_{t=1}^{n}
f^2 \bigl[(x_t + x )/h\bigr] - g(x) \int
f^2(s) \,\mrmd s \Biggr| = \RMO_P \biggl[h^r + \sqrt{
\frac{\log
n}{nh} } \biggr],
\end{equation}
which yields (\ref{12}) with $a(n)=n$ and (\ref{12a}) with $a(n)=n$
and $b_n=C_0$, where $C_0$ is a constant such that $\inf_{|x|\le C_0} g(x)>0$.
\end{rem}
%
%
\begin{rem} \label{rem25}It is much more complicated if $x_t$ is a
null recurrent Markov chain, even in the simple situation that $x_t$ is
a random walk defined as in Corollary \ref{th3}. In this regard, we
have (\ref{ez1}), but it is not clear at the moment if it is possible
to establish a result like
%
%
\begin{equation}\label{ez2}
\sup_{|x|\le b_n } \Biggl| \frac{1}{\sqrt nh} \sum_{t=1}^{n}
f^2 \bigl[(x_t + x )/h\bigr] - \int f^2(s) \,\mrmd s
L_W(1,x) \Biggr| = \RMO_P (c_n)
\end{equation}
for some $b_n\to\infty$ and $c_n\to0$. Note that (\ref{ez2})
implies that
%
%
\begin{equation}\label{ez3}
\frac{1}{\sqrt nh} \sum_{t=1}^{n}
f^2 \bigl[(x_t + y )/h\bigr]\to_P\int
f^2(s) \,\mrmd s L_W(1,0)
\end{equation}
for any fixed $y$. This is a stronger convergence than that given in
(\ref{ez1}). Our experiences show that it might not be possible to
prove (\ref{ez2}) without enlarging the probability space in which the
$x_t$ hosts.
\end{rem}

\section{Applications in non-linear cointegrating regression} \label{sec3}

Consider a non-linear cointegrating regression model:
%
%
\begin{equation}\label{mo1}
y_{t}=m(x_{t})+u_{t},\qquad t=1,2,\ldots,n,
\end{equation}
where $u_{t}$ is a stationary error process and $x_{t}$ is a
non-stationary
regressor. Let $K(x)$ be a non-negative real function and set\vadjust{\goodbreak} $%
K_{h}(s)=h^{-1}K(s/h)$ where $h\equiv h_{n}\rightarrow0$. The
conventional kernel estimate of $m(x)$ in model (\ref{mo1}) is given
by
%
%
\begin{equation}\label{add9}
\hat{m}(x)=\frac{\sum_{t=1}^{n}y_{t}K_{h}(x_{t}-x)}{%
\sum_{t=1}^{n}K_{h}(x_{t}-x)}.
\end{equation}
The point-wise limit behavior of $\hat{m}(x)$ has currently been
investigated by many authors. Among them, Karlsen, Myklebust and Tj{\o}stheim
\cite{karlsen2} discussed the situation where $x_{t}$ is a recurrent
Markov chain. Wang and Phillips \cite{wangphillips2,wangphillips3} and
Cai, Li and Park \cite{cai} considered an alternative treatment by
making use of local time limit theory and, instead of recurrent Markov
chains, worked with partial sum representations of the type
$x_t=\sum_{j=1}^t\xi_j$ where $\xi_j$ is a general linear process. In
another paper, Wang and Phillips \cite{wangphillips2}
considered the errors $u_t$ to be serially dependent and cross
correlated with the regressor $x_t$ for small lags. For other related
works, we refer to Kasparis and Phillips \cite{kasparis},
Park and Phillips \cite{park1,park2}, Gao \textit{et al.} \cite
{gao1,gao2}, Marmer \cite{marmer},
Chen, Li and Zhang \cite{chenetal}, Wang and Phillips
\cite{wangphillips4} and Wang \cite{wang}.

This section provides a uniform convergence for the $\hat{m}(x)$ by
making direct use of Theorems \ref{th4} and \ref{th1} in developing
the asymptotics. For reading convenience, we list the assumptions as follows.
%
%
\begin{assumption}\label{assumption3.1}
\textup{(i)} $\{x_k\}_{k\ge0}$ is a $\beta$-regular
Harris recurrent Markov chain defined as in Section~\ref{sec3},
where the invariant measure $\pi$ has a bounded density function
$p(s)$ on $R$;
\textup{(ii)} $\{u_t, {\mathcal F}_t\}_{t\ge1}$ is a martingale difference,
where ${\mathcal F}_t=\sigma(x_1,\ldots, x_{t+1}, u_1,\ldots,u_t)$, satisfying
$ \sup_{t\ge1}E(|u_t|^{2p}\mid{\mathcal F}_{t-1})<\infty$, where
$p\ge1+1/\varepsilon_0$ for some $0<\varepsilon_0<\beta$.
\end{assumption}
%
%
\begin{assumption}\label{assumption3.2}
The kernel $K$ satisfies that
$\int_{-\infty}^{\infty}K(s)\,\mrmd s<\infty$, $\sup_xK(x)<\infty$
and for any $x, y \in R$,
\[
\bigl|K(x)-K(y)\bigr| \le C |x-y|.
\]
\end{assumption}
%
%
\begin{assumption}\label{assumption3.3}
There exists a real positive function $g(x)$ such that
\[
\bigl|m(y)-m(x)\bigr| \leq C |y-x|^{\alpha} g(x),
\]
uniformly for some $0<\alpha\le1$ and any $(x, y)\in\Omega_\varepsilon$,
where $\varepsilon$ can be chosen sufficiently small and
$\Omega_{\varepsilon} = \{(x,y)\dvtx |y-x|\le\varepsilon, x\in R\}$.
\end{assumption}

Assumption \ref{assumption3.1} is similar to, but weaker than those
appeared in
Karlsen, Myklebust and Tj{\o}stheim \cite{karlsen2}, where the
authors considered the point-wise convergence in distribution.

Assumption \ref{assumption3.2} is a standard condition on $K(x)$ as in
the stationary
situation. The Lipschitz condition on $K(x)$ is not necessary if we
only investigate the point-wise asymptotics. See Remark~\ref{rem32} for
further details.

Assumption \ref{assumption3.3} requires a Lipschitz-type condition in a small
neighborhood of the targeted set for the functionals to be
estimated. This condition is quite weak, which may host a wide set
of functionals. Typical examples include that $m(x)=\theta_1+\theta
_2x+\cdots+\theta_kx^{k-1}$;
$m(x)=\alpha+ \beta x^{\gamma}$;
$m(x)=x(1+\theta x)^{-1}I(x\ge0)$;
$m(x)=(\alpha+\beta \RMe ^{x})/(1+\RMe^x)$.\vadjust{\goodbreak}

We have the following asymptotic results.
%
%
\begin{theorem} \label{th31} Suppose Assumptions \ref
{assumption3.1}--\ref{assumption3.3} hold, $h\to0$
and $n^{-\varepsilon_0}a(n)h\to\infty$ where $0<\varepsilon_0<\beta$ is
given as
in Assumption \ref{assumption3.1}.
It follows that
%
%
\begin{equation}\label{q1}
\sup_{|x|\le b_n'}\bigl|\hat{m}(x)-m(x)\bigr|= \RMO_{P} \bigl\{ \bigl[a(n)h
\bigr]^{-1/2} \log^{1/2}n +h^{\alpha} \delta_n
\bigr\},
\end{equation}
where $b_n'\le b_n$, $\delta_n=\sup_{|x|\le b_n'}g(x)$ and $b_n$
satisfies that
\[
\inf_{|x|\le b_n+1} \sum_{k=1}^n E K
\bigl[(x_k+x)/h\bigr]\ge a(n) h/C_0
\]
for some $C_0>0$ and all $n$ sufficiently large. In particular, for the
random walk $x_t$ defined as in Corollary \ref{th3}, we have
%
%
\begin{equation}\label{q1a}
\sup_{|x|\le b_n'}\bigl|\hat{m}(x)-m(x)\bigr|= \RMO_{P} \bigl\{
\bigl(nh^{2} \bigr)^{-1/4} \log^{1/2}n
+h^{\alpha} \delta_n \bigr\},
\end{equation}
where $b_n'\le\tau_n\sqrt n$ for some $0<\tau_n\to0$ and $\delta_n=\sup
_{|x|\le b_n'}g(x)$.
\end{theorem}
%
%
\begin{rem} \label{rem31} When a high moment exists on the error $u_t$,
the $\varepsilon_0$ can be chosen sufficiently small so that there are more
bandwidth choices in practice. It is understandable that the results
(\ref{q1})
and (\ref{q1a}) are meaningful if only $h^{\alpha}\delta_n\to0$, which
depends on the tail of the unknown regression function $m(x)$, the
bandwidth $h$ and the range $|x|\le b_n'$.
When $m(x)$ has a light tail such as $m(x)=(\alpha+\beta
\RMe^{x})/(1+\RMe^x)$, $\delta_n$ may be bounded by a constant. In this
situation, the $b_n'$ in (\ref{q1a}) can be chosen to be $\tau_n\sqrt
n$ for some $0<\tau_n\to0$.
In contrast to Theorem \ref{th1} and Remark \ref{rem23}, this kind
of range $|x|\le\tau_n\sqrt n$ might be optimal, that is, the $b_n'$
cannot be improved to $d \sqrt n$, for any $d>0$, to establish the
same rate of convergence as in (\ref{q1a}).
\end{rem}
%
%
\begin{rem} \label{rem32} Both results (\ref{q1}) and (\ref{q1a})
are sharp. However, a better result can be obtained if we are only
interested in the point-wise asymptotics for $\hat m(x)$. For instance,
as in
Wang and Phillips \cite{wangphillips1,wangphillips2} with minor
modification, we may show that, for each $x$,
%
%
\begin{equation}
\hat m(x)-m(x) = \RMO_{P} \bigl\{\bigl(nh^2
\bigr)^{-1/4}+h^{\alpha} \bigr\},
\end{equation}
whenever $x_t$ is a random walk defined as in Corollary \ref{th3}.
Furthermore $\hat m(x)$ has an asymptotic distribution that is
mixing normal, under minor additional conditions. More details are
referred to Wang and Phillips \cite{wangphillips1,wangphillips2}.
\end{rem}
%
%
\begin{rem} \label{rem33} Wang and Wang \cite{wangwang} established
a similar result to (\ref{q1a}) with the $x_t$ being a partial sum of
linear process, but only for the $x$ being a compact support and
imposing a bounded condition on $u_t$. The setting on the $x_t$ in this
paper is similar to that given in Gao, Li and Tj{\o}stheim~\cite
{gao3}, but our result provides the optimal range for the uniform
convergence holding true and removes the independence between the error
$u_t$ and $x_t$ required by Gao, Li and
Tj{\o}stheim~\cite{gao3}.\vadjust{\goodbreak}
\end{rem}

\section{Proofs of main results} \label{sec4}

\begin{pf*}{Proof of Theorem \ref{th4}} We split the set $A_n=\{x\dvtx
\|
x\|\le b_n\}$ into $m_n$ balls of the form
\[
A_{nj} =\bigl\{x\dvtx \|x-y_j\|\le1/m_n'
\bigr\},
\]
where $m_n'=[nh^{-1}/(c_n\log n)^{1/2}]$, $m_n=(b_nm_n')^d$ and $y_j$
are chosen so that $A_n\subset\bigcup A_{nj}$.
It follows that
%
%
\begin{eqnarray}\label{p6}
&&\sup_{\|x\|\le b_n} \Biggl|\sum_{t=1}^{n}
u_t f\bigl[(x_t+x)/h\bigr] \Biggr|
\nonumber
\\
&&\quad\le \max_{0\le j\le m_n}\sup_{x\in A_{nj}} \sum
_{t=1}^{n} |u_t| \bigl|f\bigl[(x_t+x)/h
\bigr]-f\bigl[(x_{t}+y_j)/h\bigr] \bigr|
\nonumber\\[-8pt]\\[-8pt]
&&\qquad{} + \max_{0\le j\le m_n} \Biggl|\sum_{t=1}^{n}
u_t f\bigl[(x_{t}+y_j)/h\bigr] \Biggr|
\nonumber
\\
&&\quad:= \lambda_{1n}+\lambda_{2n}.\nonumber
\end{eqnarray}
Recalling the Assumption \ref{assumption2.2}, it is readily seen that
%
%
\begin{eqnarray}\label{p7}
\lambda_{1n} &\le& \sum_{t=1}^n|u_t|
\max_{0\le j\le m_n}\sup_{x\in A_{nj}} \bigl|f\bigl[(x_t+x)/h\bigr]-f
\bigl[(x_{t}+y_j)/h\bigr] \bigr|
\nonumber
\\
&\le& C \bigl(h m_n'\bigr)^{-1} \sum
_{t=1}^n|u_t|
\\
&\le& C (c_n \log n)^{1/2} \frac1n \sum
_{t=1}^n|u_t| = \RMO\bigl[(c_n
\log n)^{1/2}\bigr] ,\qquad\mbox{a.s.}\nonumber
\end{eqnarray}
by the strong law of large number.

In order to investigate $\lambda_{2n}$, write $u_t'=u_tI[|u_t|\le
(c_n/\log n)^{1/2}]$ and $u_t^*=u_t'-E(u_t'\mid{\mathcal F}_{t-1})$. Recalling
$E(u_t\mid{\mathcal F}_{t-1})=0$ and $\sup_x|f(x)|<\infty$, we have
%
%
\begin{eqnarray}\label{p8}
\lambda_{2n} &\le& \max_{0\le j\le m_n} \Biggl|\sum
_{t=1}^{n} u_t^* f\bigl[(x_{t}+y_j)/h
\bigr] \Biggr|
\nonumber
\\
&&{} +\max_{0\le j\le m_n} \Biggl|\sum_{t=1}^{n}
\bigl[\bigl|u_t-u_t'\bigr|+E\bigl(\bigl|u_t-u_t'\bigr|
\mid{\mathcal F}_{t-1}\bigr) \bigr] f\bigl[(x_{t}+y_j)/h
\bigr] \Biggr|
\nonumber\\[-8pt]\\[-8pt]
&\le& \max_{0\le j\le m_n} \Biggl|\sum_{t=1}^{n}
u_t^* f\bigl[(x_{t}+y_j)/h\bigr] \Biggr| + C \sum
_{t=1}^{n} \bigl[\bigl|u_t-u_t'\bigr|+E
\bigl(\bigl|u_t-u_t'\bigr|\mid{\mathcal
F}_{t-1}\bigr) \bigr]
\nonumber
\\
:\!\!&=&\lambda_{3n} +\lambda_{4n}.\nonumber
\end{eqnarray}
Routine calculations show that, under $\sup_{t\ge1}E(|u_t|^{2p}\mid
{\mathcal F}_{t-1})<\infty$ and $n c_n^{-p} \log^{p-1}n=\RMO(1)$,
%
%
\begin{eqnarray}\label{p9}
\lambda_{4n} &\le& \sum_{t=1}^n
\bigl[ |u_t|I\bigl\{|u_t|> (c_n/\log
n)^{1/2} \bigr\} + E \bigl( |u_t| I\bigl\{|u_t|
> (c_n/ \log n )^{1/2} \bigr\} \mid\mathcal
F_{t-1} \bigr) \bigr]
\nonumber
\\
&\le& C \biggl(\frac{c_n}{\log n} \biggr)^{(1-{2p})/2} \sum
_{t=1}^n \bigl[|u_t|^{2p}+E
\bigl(|u_t|^{2p}\mid{\mathcal F}_{t-1}\bigr)
\bigr]
\nonumber\\[-8pt]\\[-8pt]
& \le& C (c_n \log n)^{1/2} \frac1n \sum
_{t=1}^n \bigl[|u_t|^{2p}+E
\bigl(|u_t|^{2p}\mid{\mathcal F}_{t-1}\bigr)
\bigr]
\nonumber
\\
&=& \RMO \bigl[(c_n \log n)^{1/2} \bigr],\qquad\mbox{a.s.}\nonumber
\end{eqnarray}
by the strong law of large number again.

We next consider $\lambda_{3n}$. As $E [(u_t^*)^2\mid{\mathcal
F}_{t-1}]\le2 (E [|u_t|^{2p}\mid{\mathcal F}_{t-1}])^{1/p}$, a.s.,
Assumptions \ref{assumption2.1} and~\ref{assumption2.3} imply that
%
%
\begin{equation}\label{*}
\max_{0\le j\le m_n} \sum_{t=1}^{n}
f^2\bigl[(x_{t}+y_j)/h\bigr]E \bigl[
\bigl(u_t^*\bigr)^2\mid{\mathcal F}_{t-1}
\bigr]=\RMO_P(c_n).
\end{equation}
Hence, for any $\eta>0$, there exists a $M_0>0$ such that
\[
P \Biggl(\max_{0\le j\le m_n} \sum_{t=1}^{n}
\sigma_{tj}^2\ge M_0 c_n\Biggr)
\le\eta,
\]
where $\sigma_{tj}^2= f^2[(x_{t}+y_j)/h]E [(u_t^*)^2\mid{\mathcal
F}_{t-1}]$, whenever $n$ is sufficiently large. This, together with
$|u_t^*|\le2(c_n /\log n)^{1/2}$ and the well-known martingale
exponential inequality (see, e.g., de la Pe{\~n}a \cite{delapena}),
implies that, for any $\eta>0$, there exists a $M_0\ge6d(k+3)$ ($k$ is
as in Assumption~\ref{assumption2.3}) such that, whenever $n$ is
sufficiently large,
%
%
\begin{eqnarray}\label{p10}
&& P\bigl[\lambda_{3n} \ge M_0 (c_n \log n
)^{1/2}\bigr]
\nonumber
\\
&&\quad\le P \Biggl[\lambda_{3n} \ge M_0 (c_n
\log n )^{1/2}, \max_{0\le j\le m_n} \sum_{t=1}^{n}
\sigma_{tj}^2\le M_0 c_n \Biggr]+
\eta
\nonumber\\[-8pt]\\[-8pt]
&&\quad\le \sum_{j=0}^{m_n} P \Biggl[\sum
_{t=1}^{n}u_t^* f
\bigl[(x_{k}+y_j)/h\bigr]\ge M_0
(c_n \log n )^{1/2}, \sum_{t=1}^{n}
\sigma_{tj}^2 \le M_0 c_n \Biggr]
+\eta
\nonumber
\\
&&\quad\le m_n \exp\biggl\{-\frac{M_0^2 c_n \log n} {
6 M_0c_n } \biggr\}+\eta\le
m_nn^{-M_0/6}+\eta\le2 \eta,\nonumber
\end{eqnarray}
where we have used the following fact:
\[
m_n \le C\bigl[n^{k+1} h^{-1}/(c_n
\log n)^{1/2}\bigr]^d \le C_1 n^{(k+2)d}
\]
as $c_n\to\infty$ and $nh\to\infty$.
This yields $\lambda_{3n}=\RMO_P [(c_n \log n)^{1/2} ]$.
Combining (\ref{p6})--(\ref{p10}), we establish (\ref{m10}).

To prove (\ref{m11}), by checking (\ref{p6})--(\ref{p9}), it
suffices to show that
%
%
\begin{equation}\label{p12}
\lambda_{3n} = \RMO \bigl[(c_n \log n)^{1/2}
\bigr],\qquad\mbox{a.s.}
\end{equation}
under the alternative condition (\ref{m1a9}). In fact, as in (\ref
{*}), it follows from (\ref{m1a9}) that
\[
\max_{0\le j\le m_n} \sum_{t=1}^{n}
f^2\bigl[(x_{t}+y_j)/h\bigr]E \bigl[
\bigl(u_t^*\bigr)^2\mid{\mathcal F}_{t-1}
\bigr]=\RMO(c_n),\qquad\mbox{a.s.}
\]
Similarly to proof of (\ref{p10}), we have for sufficiently large
$M_0$ ($M_0\ge6d (k+4)$, say),
%
%
\begin{eqnarray}\label{p13}
&& P\bigl[\lambda_{3n} \ge M_0 (c_n \log n
)^{1/2}\mbox{, i.o.}\bigr]
\nonumber
\\
&&\quad= P \Biggl[\lambda_{3n} \ge M_0 (c_n \log
n )^{1/2}, \max_{0\le j\le m_n} \sum_{k=1}^{n}
\sigma_k^2\le M_0 c_n\mbox{, i.o.} \Biggr]
\nonumber\\
&&\quad\le \lim_{s\to\infty} \sum_{n=s}^{\infty}
P \Biggl[\lambda_{3n} \ge M_0 (c_n \log n
)^{1/2}, \max_{0\le j\le m_n} \sum_{k=1}^{n}
\sigma_k^2\le M_0 c_n \Biggr]
\\
&&\quad\le\lim_{s\to\infty} \sum_{n=s}^{\infty}
m_n \exp\biggl\{ -\frac{M_0^2 c_n \log n} {
6 M_0c_n } \biggr\}
\nonumber
\\
&&\quad\le C \lim_{s\to\infty} \sum_{n=s}^{\infty}
n^{(k+2)d} n^{-M_0/6}=0,\nonumber
\end{eqnarray}
which yields (\ref{p12}). The proof of Theorem \ref{th4} is now complete.
\end{pf*}
\begin{pf*}{Proof of Theorem \ref{th4a}}
We only prove (\ref{m10a}). It is similar to prove (\ref{m11a}) and
hence the details are omitted.
We may write
%
%
\begin{eqnarray}\label{h1}
&&
\sum_{t=1}^{n} u_t f
\bigl[(x_t+x)/h\bigr]\nonumber\\
&&\quad=\sum_{t=1}^{n}
u_t f\bigl[(x_t+x)/h\bigr]I\bigl(\|x_t\|\le
b_n/2 \bigr)
\nonumber\\[-8pt]\\[-8pt]
&&\qquad{} +\sum_{t=1}^{n} u_t f
\bigl[(x_t+x)/h\bigr]I\bigl(\|x_t\|>b_n/2 \bigr)
\nonumber\\
&&\quad=\lambda_{5n} (x)+\lambda_{6n}(x)\qquad \mbox{say}.\nonumber
\end{eqnarray}
It is readily seen from (\ref{m10}) and $n \sup_{\|x\|> b_n/2}
|f(x/h)| =\RMO [(c_n \log n)^{1/2} ]$ that
\begin{eqnarray*}
\sup_{x\in R^d}\bigl|\lambda_{5n} (x)\bigr| &\le& \sup_{\|x\|\le b_n}\bigl|
\lambda_{5n} (x)\bigr| +\sup_{\|x\|> b_n}\bigl|\lambda_{5n} (x)\bigr|
\\
&\le& \RMO_P \bigl[(c_n \log n)^{1/2} \bigr]\\
&&{} +
\sup_{\|x\|> b_n/2} \bigl|f(x/h)\bigr|\sum_{t=1}^{n}
|u_t|
\\
&\le& \RMO_P \bigl[(c_n \log n)^{1/2} \bigr]
\end{eqnarray*}
as $ \frac1n\sum_{t=1}^{n} |u_t|=\RMO(1)$, a.s. by the strong law. As
for $\lambda_{6n}(x)$, we have
\begin{eqnarray*}
E \sup_{x\in R^d}\bigl|\lambda_{6n} (x)\bigr| &\le& C \sum
_{t=1}^{n} E \bigl[|u_t| I\bigl(
\|x_t\|>b_n/2 \bigr) \bigr]
\\
&\le& C \sum_{t=1}^{n} P\bigl(\|x_t
\|>b_n/2 \bigr) \le C b_n^{-k_0} \sum
_{t=1}^n E\|x_t\|^{k_0}
\\
&=& \RMO \bigl[ (c_n \log n)^{1/2} \bigr],
\end{eqnarray*}
which yield $\sup_{x \in R^d} |\lambda_{6n} (x)|=\RMO_P [ (c_n \log
n)^{1/2} ]$.
Taking these estimates into (\ref{h1}), we obtain (\ref{m10a}).
The proof of Theorem \ref{th4a} is complete.
\end{pf*}
\begin{pf*}{Proof of Theorem \ref{th1}}
First, assume there exists a $C\in\mathcal{E}^+$ such that
%
%
\begin{equation}\label{1}
P(x, A) \ge bI_C(x) \nu(A),\qquad x\in E, A\in\mathcal{E},
\end{equation}
for some $b>0$ and probability measure $\nu$ on $(E, \mathcal{E})$
with $\nu(C)>0$. Under this addition assumption, Theorem
\ref{th1} can be established by using the so-called split chain
technique. To this end, define new random variables
$Y_0, Y_1,\ldots$ and $\bar x_0, \bar x_1,\ldots$ by
\begin{eqnarray*}
P(\bar x_0\in A) &=&\nu(A),
\\
P(Y_n=1\mid\bar x_n=x) &=&h(x),
\\
P(Y_n=0\mid\bar x_n=x) &=&1-h(x),
\\
P(\bar x_{n+1}\in A\mid\bar x_n=x, Y_n=1) &=&
\nu(A),
\\
P(\bar x_{n+1}\in A\mid\bar x_n=x, Y_n=0) &=&
\frac{P(x, A)-h(x)\nu
(A)}{1-h(x)},
\end{eqnarray*}
where $h(x)=bI_C(x)$. As easily seen, $\{\bar x_n, Y_n\}_{n=0}^{\infty}$
is a Harris recurrent Markov chain with state space $E\times\{0,1\}$
and $\{\bar x_n\}_{n=0}^{\infty}$ has
the same transition probability $P(x, A)$ as those of $\{x_n\}
_{n=0}^{\infty}$. Since our result is free of the initial distribution,
$\{x_n\}_{n=0}^{\infty}$ can be assumed to be identical with $\{\bar
x_n\}_{n=0}^{\infty}$, that is, $x_0$ has the distribution $\nu$.

Further define $\rho_0=-1$,
\[
\rho_k=\min\{i\dvtx i\ge\rho_{k-1}, Y_i=1\},\qquad
k=1,2,\ldots,
\]
$N(n)=\max\{k\dvtx \rho_k\le n\}$, and
\[
Z_{j}(x) = \sum_{k=\rho_{j-1}+1}^{\rho_{j}}
f^2\bigl[(x_k+x)/h\bigr],\qquad Z_{jn}(x) = \sum
_{k=\rho_{j-1}\wedge n+1}^{\rho_{j}\wedge
n} f^2
\bigl[(x_k+x)/h\bigr]
\]
for $j=1, 2,\ldots\,$. It is well known that the blocks
\[
(x_{\rho_i+1},\ldots,x_{\rho_{i+1}}),\qquad i=0, 1,2,\ldots,
\]
are i.i.d. blocks, $x_{\rho_i+1}$ having the distribution $\nu$.
Hence, for each $h$ and $x$, $\{Z_j^*(x), \rho_j-\rho_{j-1}\}
_{j=1}^{\infty}$, where $Z_j^*(x)=Z_j(x)$ or $Z_{jn}(x)$ is a sequence
of i.i.d. random vectors. Furthermore, by recalling that $\pi$ has a
bounded density function $p(s)$, $\int_{-\infty}^{\infty}|f(x)|\,\mrmd \pi
(x)<\infty$ and $\sup_s|f(s)|<\infty$, we have
%
%
\begin{eqnarray}\label{45}
EZ_1(x) &=& b \int_{-\infty}^{\infty}
f^2\bigl[(s+x)/h\bigr] \pi(\mrmd s)
\nonumber\\[-8pt]\\[-8pt]
&=&b h \int_{-\infty}^{\infty} f^2(s) p(-x+sh)
\,\mrmd s \le C^* h\nonumber
\end{eqnarray}
for any $x\in R$ and
%
%
\begin{equation}\label{fi2}
\sup_{x\in R}E\bigl|Z_1(x)\bigr|^{2k} \le C h
\end{equation}
for any integer $k$.
See Lemma 5.2 of Karlsen and Tj{\o}stheim \cite{karlsen1} or Lemma
B.1 of Gao, Li and Tj{\o}stheim \cite{gao3}.
We also have the following lemma.
%
%
\begin{lem} \label{lem1}Suppose that $d_n\sim C_0a(n)$, where $C_0>0$
is a constant, and all $y_j, j=0,1,\ldots,m_n$, are different, where
$|y_j|\le n^{m_0}$ and $m_n\le n^{m_1}$ for some $m_0, m_1>0$. Then,
%
%
\begin{eqnarray}
\label{95}
R_{n} &:=& \max_{0\le j\le m_n} \Biggl|\sum_{k=0}^{d_n}
\bigl[Z_k^*(y_j)-EZ_k^*(y_j)
\bigr] \Biggr| =\RMO_P\bigl[n^{-\varepsilon_0/4} a(n) h\bigr],
\\
\label{96}
\Delta_{n} &:=& \max_{0\le j\le m_n}E \Biggl|\sum
_{k=0}^{d_n} \bigl[Z_k^*(y_j)-EZ_k^*(y_j)
\bigr] \Biggr| =\RMO\bigl[n^{-\varepsilon_0/4} a(n) h\bigr],
\end{eqnarray}
where $\varepsilon_0$ is a constant such that $n^{-\varepsilon
_0}a(n)h\to
\infty$.
\end{lem}
\begin{pf} Only\vspace*{1pt} consider $Z_k^*(x)=Z_k(x)$, as the situation that
$Z_k^*(x)=Z_{kn}(x)$ is similar.
To this end, write $\widetilde Z_{ i}(y_j)= Z_{ i}(y_j)I(|Z_{
i}(y_j)|\le n^{-\varepsilon_0/2} a(n)h)$ and
$\widehat Z_{ i}(y_j)= Z_{ i}(y_j)I(|Z_{ i}(y_j)|> n^{-\varepsilon_0/2}
a(n) h)$. We have
%
%
\begin{eqnarray}\label{b1}
R_n & \le& \max_{0\le j\le m_n} \Biggl|\sum_{i=1}^{d_n}
\bigl[\widetilde Z_{ i}(y_j)-E\widetilde
Z_{ i}(y_j) \bigr] \Biggr| + \max_{0\le j\le
m_n} \sum
_{i=1}^{d_n} \bigl[ \widehat Z_{ i}(y_j)+
E \widehat Z_{
i}(y_j) \bigr]
\nonumber\\[-8pt]\\[-8pt]
:\!\!&= & R_{1n}+R_{2n}.\nonumber
\end{eqnarray}
By taking $k\ge(m_1+2)/\varepsilon_0$ in (\ref{fi2}) and noting
$n^{-\varepsilon
_0}a(n)h\to\infty$, simple calculations show that
%
%
\begingroup
\abovedisplayskip=7pt
\belowdisplayskip=7pt
\begin{eqnarray}\label{i9}
E R_{2n} &\le& C m_n a(n)\max_{0\le j\le m_n}E
Z_{ 1}(y_j)I\bigl(\bigl|Z_{
i}(y_j)\bigr|>
n^{-\varepsilon_0/2}a(n)h\bigr)
\nonumber
\\
&\le& C_{1} a(n) h \bigl( n^{m_1+1-k\varepsilon_0}h^{-1} \bigr)\le
C_{1} a(n) h (nh)^{-1}
\\
&\le& C_{2} n^{-\varepsilon_0/2}a(n) h,\nonumber
\end{eqnarray}
which yields $R_{2n}=\RMO_P[n^{-\varepsilon_0/2} a(n) h]$. As for
$R_{1n}$, by
using (\ref{fi2}) with $k=2$ and noting
\[
E\RMe^{t(\widetilde Z_{ i}(y_j)-E\widetilde Z_{ i}(y_j))} \le1+ \frac
{t^2}2 EZ_{ 1}^2(y_j)\RMe^{2tn^{-\varepsilon_0/2}a(n)h}
\le \RMe^{C_0 t^2 h}
\]
for any $t\le(n^{-\varepsilon_0/2} a(n)h)^{-1}$ and some $C_0>0$, the
standard Markov inequality implies that
%
%
\begin{eqnarray}\label{l1}
&&P\bigl(R_{1n} \ge M n^{-\varepsilon_0/4} a(n)h\bigr)
\nonumber
\\
&&\quad\le C m_n \max_{0\le j\le m_n} P \Biggl( \Biggl|\sum
_{i=1}^{C_{\varepsilon} a(n)} \bigl[\widetilde Z_{ i}(y_j)-E
\widetilde Z_{
i}(y_j) \bigr] \Biggr|\ge M n^{-\varepsilon_0/4}a(n)h
\Biggr)
\nonumber\\[-8pt]\\[-8pt]
&&\quad\le C m_n \exp\bigl(-M t n^{-\varepsilon_0/4} a(n)h
+C_{\varepsilon
}a(n)t^2h\bigr)
\nonumber
\\
&&\quad\le C m_n \exp\bigl(-M n^{\varepsilon_0/4}/4 \bigr)\to0\nonumber
\end{eqnarray}
as $n\to\infty$. Hence, $R_{1n}=\RMO_P[n^{-\varepsilon_0/4} a(n) h]$.
Combining (\ref{b1})--(\ref{l1}), we prove (\ref{95}).

The proof of (\ref{96}) is similar except more simpler. Indeed, by
independence of $\widetilde Z_{ i}(x)$, we obtain
\begin{eqnarray*}
\Delta_n & \le& \max_{0\le j\le m_n} E \Biggl|\sum
_{i=1}^{d_n} \bigl[\widetilde Z_{ i}(y_j)-E
\widetilde Z_{ i}(y_j) \bigr] \Biggr| \\
&&{}+ 2\max_{0\le j\le
m_n}
\sum_{i=1}^{d_n}E \widehat
Z_{ i}(y_j)
\\
& \le& 2 \max_{0\le j\le m_n} d_n^{1/2} \bigl[E\widetilde
Z_{
1}^2(y_j)\bigr]^{1/2} + C
n^{-\varepsilon_0/2}a(n) h
\\
&\le& C n^{-\varepsilon_0/4}a(n) h,
\end{eqnarray*}
due to the fact:
\[
E\widetilde Z_{ 1}^2(y_j) \le
n^{-\varepsilon_0/2} a(n)h E Z_1(y_j) \le C
n^{-\varepsilon_0/2} a(n)h^2.
\]
The proof of Lemma \ref{lem1} is complete.
\end{pf}\eject
\endgroup

We are now ready to prove (\ref{12}) and (\ref{12a}) under the
additional condition (\ref{1}).

(\ref{12}) first. As in proof of (\ref{p6}) and (\ref{p7}), but
letting $ y_j=-[n^m]-1+j / m_n',j=0, 1,2,\ldots, m_n, $ where
$m_n'=[nh^{-2}/a(n)]$ and $m_n=2([n^m]+1)m_n'$, we have
%
%
\begin{equation}\label{fi1}
\sup_{|x|\le n^m} \sum_{k=0}^n
f^2\bigl[(x_k+x)/h\bigr] \le\max_{0\le j\le
m_n}\sum
_{k=0}^n f^2
\bigl[(x_k+y_j)/h\bigr]+ C a(n) h.
\end{equation}
Note that
\[
\sum_{k=0}^n f^2
\bigl[(x_k+x)/h\bigr] \le\sum_{k=0}^{\rho_{ N(n+1)}}
f^2\bigl[(x_k+x)/h\bigr] = \sum
_{i=1}^{N(n+1)} Z_{ i}(x),
\]
and $ \{N(n)/a(n) \}_{n\ge1}$ is bounded in probability. See,
for example, Chen \cite{chen}. For each $\varepsilon>0$, there exist
$0<C_{\varepsilon},
C_{1\varepsilon}<\infty$ such that
%
%
\begin{equation}\label{a99}
P\bigl(C_{1\varepsilon}a(n)\le N(n)\le C_{\varepsilon}a(n)\bigr) \ge1-
\varepsilon,
\end{equation}
whenever $n$ is sufficiently large. Consequently, for each $a>0,
\varepsilon
>0$ and $n$ large enough,
\[
P \Biggl(\max_{0\le j\le m_n}\sum_{k=0}^n
f^2\bigl[(x_k+y_j)/h\bigr]\ge a \Biggr)\le P
\Biggl(\max_{0\le j\le m_n} \sum_{i=1}^{C_\varepsilon a(n)}
Z_{
i}(y_j)\ge a \Biggr)+\varepsilon.
\]
This, together with (\ref{95}) with $Z_k^*(x)=Z_k(x)$, implies (\ref
{12}) under (\ref{1}), since
\[
\max_{0\le j\le m_n} \sum_{i=1}^{C_\varepsilon a(n)}
Z_{ i}(y_j) \le C_\varepsilon a(n)
\max_{0\le j\le m_n}E Z_{ 1}(y_j)+R_n =
\RMO_P\bigl[a(n) h\bigr].
\]

We next consider (\ref{12a}) under (\ref{1}). To this regard,
let $ y_j=-[b_n]-1+j / m_n', j=0, 1,2,\ldots, m_n,
$ where $m_n'=[ n^{1+\varepsilon_0/2}h^{-2}/a(n)]$ and $m_n=2([b_n]+1)m_n'$.
Since
\begin{eqnarray*}
&& \max_{0\le j\le m_n-1}\sup_{x\in
[y_j, y_{j+1}]} \sum_{t=1}^{n}
\bigl|f^2\bigl[( x_t+x)/h\bigr]-f^2\bigl[(
x_t+y_j)/h\bigr] \bigr|
\\
&&\quad \le C n h^{-1} \max_{0\le j\le
m_n-1}|y_{j+1}-y_j|
\le Cn^{-\varepsilon_0/2} a(n)h,
\end{eqnarray*}
it is readily seen that
%
%
\begin{equation}\label{145}
\inf_{|x|\le b_n} \sum_{t=1}^{n}
f^2\bigl[(x_{t}+x)/h\bigr] \ge\Delta_{1n}-\RMO_{P}
\bigl[n^{-\varepsilon_0/2} a(n)h\bigr],
\end{equation}
where $\Delta_{1n}=\inf_{1\le j\le m_n} \sum_{t=1}^{n} f^2[( x_t+y_j)/h]$.
Furthermore, by recalling (\ref{a99}) and noting that
\[
\sum_{k=0}^n f^2
\bigl[(x_k+x)/h\bigr] \ge\sum_{k=0}^{\rho_{ N(n)}\wedge n}
f^2\bigl[(x_k+x)/h\bigr] = \sum
_{i=1}^{N(n)} Z_{ in}(x),
\]
we have, for each $a>0, \varepsilon>0$ and $n$ large enough,
%
%
\begin{equation}\label{146}
P (\Delta_{1n}\ge a )\ge P \Biggl(\inf_{0\le j\le m_n} \sum
_{i=1}^{C_{1\varepsilon} a(n)} Z_{
in}(y_j)\ge a
\Biggr)-\varepsilon.
\end{equation}
On the other hand, it follows from (\ref{95}) with
$Z_k^*(x)=Z_{kn}(x)$ that
%
%
\begin{eqnarray}\label{147}
\inf_{0\le j\le m_n} \sum_{i=1}^{C_{1\varepsilon} a(n)}
Z_{ in}(y_j) &\ge& C_{1\varepsilon} a(n)
\inf_{0\le j\le m_n}EZ_{ 1n}(y_j) -R_n
\nonumber\\[-8pt]\\[-8pt]
&\ge& C_{1\varepsilon} a(n)\inf_{0\le j\le m_n}EZ_{ 1n}(y_j)
-\RMO_{P}\bigl[n^{-\varepsilon_0/4} a(n)h\bigr].\nonumber
\end{eqnarray}
Combining (\ref{145})--(\ref{147}), the result (\ref{12a}) under
(\ref{1}) will follow if we prove: there exists a $b_0>0$ such that
%
%
\begin{equation}\label{148}
\inf_{0\le j\le m_n}EZ_{ 1n}(y_j) \ge b_0 h
\end{equation}
for all $n$ sufficiently large. To prove (\ref{148}), first note that
there exists a $b_1>0$ such that $EN^2(n)/a^2(n)\le b_1$. See Lemma 3.3
of Karlsen and Tj{\o}stheim \cite{karlsen1}, for instance. Therefore,
by taking $d_n=[b_2a(n)]+1$, where $b_2>b_1$ is chosen later, we have
for some $b_0>0$,
\begin{eqnarray*}
\inf_{0\le j\le m_n} E Z_{ 1n}(y_j) &=&
\frac{1}{d_n} \inf_{0\le
j\le m_n} E\sum_{i=1}^{d_n}
Z_{ in}(y_j) \\
&=& \frac{1}{d_n} \inf_{0\le
j\le m_n} E
\sum_{t=1}^{\rho_{d_n}\wedge n} f^2\bigl[(
x_t+y_j)/h\bigr]
\\
&\ge&\frac{1}{d_n} \inf_{0\le j\le m_n}E \Biggl(\sum
_{t=1}^{ n} f^2\bigl[(
x_t+y_j)/h\bigr] -I(\rho_{d_n} \le n) \sum
_{t=1}^{ \rho_{d_n}} f^2\bigl[(
x_t+y_j)/h\bigr] \Biggr)
\\
&\ge&\frac{1}{d_n} \Biggl(\inf_{|x|\le b_n+1}E\sum
_{k=1}^nf^2\bigl[( x_t+x)/h
\bigr]-M_n \Biggr)
\\
&\ge&\frac{1}{d_n} \bigl[a(n)h/C_0-M_n \bigr]
\\
&\ge& b_0 h,
\end{eqnarray*}
whenever $n$ is sufficiently large, where we have used the condition
(\ref{wond}) and the fact: it follows from
(\ref{45}), (\ref{96}) and $\rho_{d_n}\le n$ if and only if $N(n)>
d_n$ that
\begin{eqnarray*}
M_n :\!\!&=& \max_{0\le j\le m_n}E \Biggl[I(\rho_{d_n} \le n)
\sum_{t=1}^{ \rho_{d_n}} f^2\bigl[(
x_t+y_j)/h\bigr] \Biggr]
\\
&= & \max_{0\le j\le m_n} E \Biggl[I(\rho_{d_n} \le n) \sum
_{i=1}^{d_n} Z_i(y_j)
\Biggr]
\\
&\le& d_n \max_{0\le j\le m_n} EZ_1(y_j) P
\bigl(N(n)\ge d_n\bigr)+ \max_{0\le j\le m_n} E \Biggl|\sum
_{i=1}^{d_n} \bigl[Z_i(y_j)-
EZ_i(y_j) \bigr] \Biggr|
\\
&\le& b_2^{-1} C^*h a^{-1}(n)
EN^2(n)+\RMO\bigl[n^{-\varepsilon_0/4}a(n)h\bigr]
\\
&\le& C_0^{-1} a(n)h/2
\end{eqnarray*}
by choosing $b_2=3 C_0 b_1 C^*$ and $n$ sufficiently large.
This proves (\ref{148}) and also completes the proof of (\ref{12a})
under (\ref{1}).

We now consider general situation. Let $0<t<1$ be fixed. Define a
transition probability $P_t(x, A)$ on
$(E, \mathcal{E})$ by
\[
P_t(x, A) =(1-t) \sum_{k=1}^{\infty}t^{k-1}
P^k(x, A),\qquad x\in E, A\in\mathcal{E}.
\]
Let $\{\beta_n\}_{n\ge1}$ be an i.i.d. Bernoulli random variables
with the common law
\[
P(\beta_1=0)=t \quad\mbox{and}\quad P(\beta_1=1)=1-t
\]
and assume $\{\beta_n\}_{n\ge1}$ and $\{x_n\}_{n\ge0}$ are
independent. Define a renewal sequence $\{\sigma(k)\}_{k\ge0}$ by
\[
\sigma(0)=0 \quad\mbox{and}\quad \sigma(k)=\inf\bigl\{n\dvtx n\ge\sigma(k-1);
\beta_n=1\bigr\},\qquad k\ge1.
\]
With these notations, $\{x_{\sigma(n)}\}_{n\ge0}$ is a Harris recurrent
Markov chain with the invariant measure~$\pi$.
The transition probability $P_t(x, A)$ of $\{x_{\sigma(n)}\}_{n\ge0}$
satisfies the additional condition (\ref{1})
and
\[
a_t(n):=\pi(D)^{-1} \sum_{k=1}^n
\nu P_t^k(D) \sim(1-t)^{1-\gamma
} a(n).
\]
See Chen \cite{chen}, for instance. By virtue of these facts, it
follows from the first part proof of (\ref{12}) that, for any fixed
$m>0$ and $h>0$,
\[
\sup_{|x|\le n^m} \sum_{k=1}^{\sigma(n)}
\beta_kf^2\bigl[(x_k+x)/h\bigr] =
\sup_{|x|\le n^m} \sum_{k=1}^n
f^2\bigl[(x_{\sigma(k)}+x)/h\bigr]= \RMO_P
\bigl[a_t(n) h\bigr].
\]
Now by noting $\sigma([\lambda n])/n\to_{\mathrm{a.s.}} \lambda/(1-t)$
by the strong
law and taking $\lambda$ such that $\lambda/(1-t)\ge1$, simple
calculations show that
\[
\sup_{|x|\le n^m} \sum_{k=1}^{n}
\beta_kf^2\bigl[(x_k+x)/h\bigr] \le
\sup_{|x|\le n^m} \sum_{k=1}^{\sigma([\lambda n])}
f^2\bigl[(x_{\sigma
(k)}+x)/h\bigr]=\RMO_P\bigl[a(n) h
\bigr].
\]
Similarly,
\[
\sup_{|x|\le n^m} \sum_{k=1}^{n} (1-
\beta_k)f^2\bigl[(x_k+x)/h\bigr]
=\RMO_P\bigl[a(n) h\bigr]
\]
and hence the result (\ref{12}) under general situation follows.

The proof of (\ref{12a}) under general situation is similar and hence
the details are omitted.
\end{pf*}
\begin{pf*}{Proof of Corollary \ref{th3}}
We first notice that:
\begin{itemize}[(F)]
\item[(F)] $x_k=\sum_{j=1}^k\varepsilon_{j}$ is a Harris null recurrent
Markov chain, satisfying (\ref{1}), $a(t)=\sqrt t$ and the invariant
measure $\pi$ is the Lebesgue measure.
\end{itemize}
Due to the fact (F), (\ref{12a1}) follows immediately from
Theorem \ref{th1}.

To prove (\ref{m15}), by Theorem \ref{th1}, it suffices to show that
(\ref{wond}) holds true with $b_n=\tau_n \sqrt n$ and $a(n)=\sqrt n$.
In fact, under the conditions of Corollary \ref{th3}, $x_k/\sqrt k$
has a density $p_k(x)$, satisfying $\sup_x|p_k(x)-\phi(x)|\to0$, as
$k\to\infty$, where $\phi(x)=\RMe^{-x^2/2}/\sqrt{2\pi}$, due to the
central limit theorem. This implies that
\[
\inf_{|x|\le3\tau_n }p_k(x)\ge\inf_{|x|\le3\tau_n }\phi(x)-
\sup_x\bigl|p_k(x)-\phi(x)\bigr| \ge A_0>0
\]
for some $A_0>0$ and all sufficiently large $k$. Hence, for $n/2<k\le
n$ and $n$ sufficiently large, we have
\begin{eqnarray*}
\inf_{|x|\le\tau_n \sqrt n+1} E f^2\bigl[(x_k+x)/h\bigr] &=&
\inf_{|x|\le
\tau_n \sqrt n+1} \int_{-\infty}^{\infty} f^2
\bigl[(\sqrt k y+x)/h\bigr] p_k(y)\,\mrmd y
\\
&\ge&\frac{h}{\sqrt k}\inf_{|x|\le\tau_n \sqrt n+1}\int_{-\infty
}^{\infty}
f^2(y) p_k\bigl[(yh-x)/\sqrt k\bigr]\,\mrmd y
\\
&\ge&\frac{h}{\sqrt k}\inf_{|x|\le3 \tau_n } p_k(x)\int
_{|y|\le
M_1} f^2(y)\,\mrmd y
\\
&\ge&\frac{A_0 h}{2\sqrt n} \int_{|y|\le M_1} f^2(y)\,\mrmd y,
\end{eqnarray*}
where $M_1$ is chosen such that $\int_{|y|\le M_1} f^2(y)\,\mrmd y>0$.
Consequently, there exists a constant $C_0>0$ such that
\[
\inf_{|x|\le\tau_n \sqrt n+1} \sum_{k=1}^n E
f^2\bigl[(x_k+x)/h\bigr]\ge\inf_{|x|\le\tau_n \sqrt n+1} \sum
_{k=n/2}^n E f^2
\bigl[(x_k+x)/h\bigr] \ge\sqrt n h/C_0
\]
as required. The proof of Corollary \ref{th3} is now complete.
\end{pf*}
\begin{pf*}{Proof of Theorem \ref{th31}}
We may write $\hat{m}(x)-m(x)$
as
%
%
\begin{eqnarray}\label{q3}
\hat{m}(x)-m(x) &=& \frac{\sum_{t=1}^{n}u_{t}K_{h}(x_{t}-x)}{%
\sum_{t=1}^{n}K_{h}(x_{t}-x)}+\frac{\sum_{t=1}^{n}
[m(x_{t})-m(x) ]%
K_{h}(x_{t}-x)}{\sum_{t=1}^{n}K_{h}(x_{t}-x)}
\nonumber\\[-8pt]\\[-8pt]
:\!\!&=& \Theta_{1n}(x)+\Theta_{2n}(x).\nonumber
\end{eqnarray}
Note that, for any $|x|\le b_n'$, there exists a $C_0>0$ such that
$K[(x_t-x)/h]=0$ if $|x_t-x|\ge hC_0$. It follows from
Assumption \ref{assumption3.3} that, whenever $n$ is sufficiently large,
\[
\sup_{|x|\le b_n'}\bigl|\Theta_{2n}(x)\bigr| \le C_1
\delta_n \sup_{|x|\le
b_n'}\frac{\sum_{t=1}^{n} |x_{t}-x|^{\alpha} K[(x_t-x)/h]}{\sum
_{t=1}^{n}K[(x_t-x)/h]} \le Ch^{\alpha}
\delta_n.
\]
This, together with (\ref{12a}) [taking $f^2(s)=K(s)$] in Theorem
\ref{th1}, implies that (\ref{q1}) will follow if we prove
%
%
\begin{equation}\label{mp1}
\sup_{|x|\le b_n} \sum_{t=1}^{n}u_{t}K
\bigl[(x_{t}-x)/h\bigr] =\RMO_P \bigl[\bigl[a(n)h
\bigr]^{1/2} \log^{1/2}n \bigr].
\end{equation}
In fact, with $p\ge1+1/\varepsilon_0$ and $c_n=a(n)h\to\infty$, we have
\[
n c_n^{-p} \log^{p-1}n\le\bigl(n^{-\varepsilon_0}a(n)h
\bigr)^{-1-1/\varepsilon
_0} n^{-\varepsilon
_0}\log^{p-1}n \to0,
\]
since $n^{-\varepsilon_0} a(n)h\to\infty$. Now, by recalling (\ref{12}),
it is readily seen that the conditions of Theorem~\ref{th4} hold for
$f(x)=K(x)$ and $c_n=a(n)h$. The result (\ref{mp1}) follows from (\ref
{m10}) in Theorem~\ref{th4}.
\end{pf*}

\section*{Acknowledgements}

The authors thank Associate Editor, two referees and Professor Jiti
Gao for helpful comments on the original version. Wang acknowledges the
partial research support from the Australian research council.



\printhistory

\end{document}